\documentclass[11pt]{article}
\usepackage{amsmath}
\usepackage[english]{babel}
\usepackage{amsfonts}
\usepackage{color}
\newtheorem{theo}{Theorem}[section]

\newtheorem{rem}[theo]{Remark}
\newtheorem{prop}[theo]{Proposition}

\newcommand{\mysection}[1]{\section{#1} \setcounter{equation}{0}}
\newcommand{\proof}{{\sc Proof.} \ }

\newcommand{\be}{\begin{equation} \label}
\newcommand{\ee}{\end{equation}}
\newcommand{\bea}{\begin{eqnarray}\label}
\newcommand{\eea}{\end{eqnarray}}
\newcommand{\bas}{\begin{eqnarray*}}
\newcommand{\eas}{\end{eqnarray*}}
\newcommand{\bit}{\begin{itemize}}
\newcommand{\eit}{\end{itemize}}
\newcommand{\qed}{\hfill$\Box$ \vskip.2cm}
\newcommand{\R}{\mathbb{R}}

\newcommand{\pO}{\partial\Omega}

\newcommand{\eps}{\varepsilon}

\newcommand{\abs}{\\[5pt]}

\newcommand{\io}{\int_\Omega}

\begin{document}
\title{Blow-up profiles for the parabolic-elliptic \\ Keller-Segel system
in dimensions $n\ge 3$}
\author
{
Philippe Souplet\footnote{souplet@math.univ-paris13.fr}\\
{\small {Universit\'{e} Paris 13, Sorbonne Paris Cit\'{e}, CNRS UMR 7539}}\\
{\small {Laboratoire Analyse, G\'{e}om\'{e}trie et Applications, 93430, Villetaneuse, France.} }
 \and
Michael Winkler\footnote{michael.winkler@math.uni-paderborn.de}\\
{\small Institut f\"ur Mathematik, Universit\"at Paderborn,}\\
{\small 33098 Paderborn, Germany} }
\date{}
\maketitle
\begin{abstract}
\noindent
We study the blow-up asymptotics of radially decreasing solutions of the parabolic-elliptic Keller-Segel-Patlak system
in space dimensions $n\ge 3$.
In view of the biological background of this system and of its mass conservation property,
blowup is usually interpreted as a phenomenon of concentration or aggregation of the bacterial population.
Understanding the asymptotic behavior of solutions at the blowup time is thus meaningful 
for the interpretation of the model.

Under mild assumptions on the initial data, for $n\ge 3$, we show that the final profile satisfies $C_1|x|^{-2}\le u(x,T)\le C_2|x|^{-2}$,
 with convergence in $L^1$ as $t\to T$.
This is in sharp contrast with the two-dimensional case, where solutions are known to 
concentrate to a Dirac mass at the origin (plus an integrable part).
We also obtain refined space-time estimates of the form $u(x,t)\le C(T-t+|x|^2)^{-1}$ for type~I blowup solutions.
Previous work had shown that radial, self-similar blowup solutions (which satisfy the above estimates)
exist in dimensions $n\ge 3$ and do not exist in dimension $2$.
Our results thus reveal that the final profile displayed by these special solutions
actually corresponds to a much more general phenomenon.
\abs
\end{abstract}
\newpage
\mysection{Introduction and main results}
Let $\Omega=B_R\subset\R^n$ or $\Omega=\R^n$, with $R>0$ and $n\ge 1$.
In this article we consider radially symmetric solutions of the well-known parabolic-elliptic Keller-Segel-Patlak system
\be{0}
        \left\{ \begin{array}{lcll}
    	\hfill u_t &=& \Delta u - \nabla \cdot (u\nabla v), 
    	& x\in\Omega, \ t>0, \\[1mm]
    	\hfill 0 &=& \Delta v + u - \mu,
    	& x\in\Omega, \ t>0, \\[1mm]
    	\hfill \frac{\partial u}{\partial\nu}&=& \frac{\partial v}{\partial\nu}\ =\ 0,
    	& x\in\pO, \ t>0, \\[1mm]
    	\hfill u(x,0)&=&u_0(x),
    	& x\in\Omega,
        \end{array} \right.
\ee
where
\be{defmu}
\mu:=
        \left\{ \begin{array}{lcll}
\displaystyle\frac{1}{|\Omega|} \io u_0&&\hbox{ if $\Omega=B_R$,} \\[3mm]
0&&\hbox{ if $\Omega=\R^n$}
        \end{array} \right.
\ee
(and the boundary conditions in (\ref{0}) are understood to be empty if $\Omega=\R^n$).
As for the initial data, unless explicitly stated otherwise, we shall always assume that
\be{i0}
\begin{array}{l}	
u_0\in L^\infty(\Omega),\ u_0\ge 0, 
\  \mbox{$u_0$ is radially symmetric} \\ [2mm]
\mbox{and nonincreasing with respect to $|x|$, with $u_0\not\equiv const$}.
\end{array}
\ee
This system (with or without radial symmetry assumptions) 
arises as a simplified model of chemotaxis, 
where $u$ and $v$ respectively stand for the density of the bacterial population 
and of the secreted chemoattractant.
System (\ref{0}) is also involved in a model of gravitational interaction of particles.
It has received considerable attention from the mathematical point of view
(see e.g. the surveys \cite{Horst1}, \cite{Horst2} and the book \cite{SuzBook} for references).

Problem (\ref{0}) is locally well posed (see Proposition~\ref{prop12} for a precise statement) and we denote by $(u,v)$ its unique, maximal classical solution,
and by $T$ its maximal existence time. 
We note that $u$ is nonnegative and radial decreasing.
It is known that if $n\ge 2$ and $u_0$ is suitably large, then $T<\infty$ and the solution blows up in the following sense:
$$	\lim_{t\to T} \|u(\cdot,t)\|_{L^\infty(\Omega)} = \infty$$ 
 (cf.~\cite{jaeger_luckhaus}, \cite{BHN}, \cite{biler1995}, 
 \cite{nagai1995}, \cite{CPZ}, \cite{Sen05}, \cite{BDP}).
For instance this is the case whenever $\|u_0\|_1>c(n)R^{n-2}$ 
 if~$\Omega=B_R$ (see Proposition~\ref{BUcond}).

Our concern in this article is about the asymptotic behavior of $u$ at  and near the blowup time.
Recalling the well-known mass conservation property $\|u(t)\|_1=\|u_0\|_1$,
and keeping in mind the biological background of system (\ref{0}),
blowup is usually interpreted as a phenomenon of concentration or aggregation of the bacterial population.
Understanding the asymptotic behavior of $u$ at the blowup time is thus meaningful 
for the interpretation of the model. 
 We focus on the properties of $u$, as the quantity $v$ can be recovered
from $u$ by the linear second equation in (\ref{0}).

For $n=2$, this behavior is rather well understood: the solution exhibits a
quantized concentration of mass at the origin. More precisely, for any blowing up solution, there exist $m\ge 8\pi$ 
and $f\in L^1(\Omega)$ such that
\be{DiracCv}
u(\cdot,t) \rightharpoonup m\delta_0+f\quad\hbox{ in ${\cal M}(\Omega)$,\quad as $t\to T$}
\ee
(see \cite{HV}, \cite{SenSu}, \cite{SenSu2}, \cite{SuzBook}).
For $n\ge 3$, the situation may be quite different.
Indeed, in the case $\Omega=\R^n$, there exist radial, positive, backward self-similar solutions of the form 
\be{Selfsimil}
u(x,t)=(T-t)^{-1}V(x/\sqrt{T-t}),\quad x\in\R^n,\ 0<t<T,
\ee
where the radially decreasing profile function $V$ satisfies $\lim_{y\to\infty}y^2V(y)=L\in (0,\infty)$
(see \cite{HMV}, \cite{Sen05}, \cite{GMS}).
This leads to the final blowup profile
\be{Homogprofile}
 U(x):=\lim_{t\to T}u(x,t)=L|x|^{-2},\quad x\ne 0,
 \ee
where the convergence also takes place in $L^1_{loc}(\R^n)$.
Therefore, these self-similar solutions have an integrable singularity and exhibit no mass concentration.
More generally, for $\Omega=B_R$ and $n\ge 2$, it is known (see, e.g., \cite{GP}) that under assumption (\ref{i0}), 
blowup can occur only at $x=0$.
However the final blowup profile of general radial decreasing solutions does not seem to have been studied so far.

\medskip
The purpose of this article is to show that the $|x|^{-2}$ final profile,
displayed by the special self-similar solutions in (\ref{Selfsimil}), actually corresponds to a much more general phenomenon,
which marks a strong difference between dimensions $n\ge 3$ and $n=2$.
In what follows, the blowup set $B(u_0)$ of $(u,v)$ is defined by 
\be{defBu0}
B(u_0):=\bigl\{x_0\in \overline\Omega;\, u(x_j,t_j)+|\nabla v(x_j,t_j)|\to \infty
\ \hbox{for some sequence $(x_j,t_j)\to (x_0,T)$}\bigr\}.
 \ee
Our first main result is the following:
\begin{theo}\label{theo13}
Let $n\ge 3$. Consider problem (\ref{0}), where $u_0$ satisfies (\ref{i0}) and $T<\infty$.
 \abs
  (i)\ Let $\Omega=B_R\subset\R^n$. Then there exists $C>0$ such that
  \be{13.1}
	u(x,t) \le \frac{C}{|x|^2},
	\qquad 0<|x|\le R,\ 0<t<T.
  \ee
Moreover, we have $B(u_0)=\{0\}$, the final blowup profile $U(x):=\lim_{t\to T}u(x,t)$ exists for all $x\in \overline\Omega\setminus\{0\}$,
 where convergence also takes place in $L^1(B_R)$, and $U$ satisfies 
  \be{13.2}
	U(x) \le \frac{C}{|x|^2}, \qquad 0<|x|\le R.
\ee
(ii)\ Let $\Omega=\R^n$ and assume $B(u_0)\neq\R^n$. 
Then assertion (i) remains valid, with $R=1$ in (\ref{13.1}) and~(\ref{13.2}).  
The assumption $B(u_0)\neq\R^n$ is satisfied in particular whenever $u_0\in L^1(\R^n)$.

\end{theo}

Our second main result shows that (\ref{13.2}) is sharp. Namely, the final blowup profile satisfies the corresponding lower estimate,
assuming the following additional hypothesis on $u_0$:
  \be{i2}
	u_0\in C^1(\overline\Omega),\ \quad r^{n-1} u_{0,r}(r) + u_0(r) \int_0^r  \big( u_0(s)-\mu\big) s^{n-1} ds \ge 0,
	\ \mbox{ for all } r\in (0,R).
\ee
We note that, for $\Omega=B_R$ and any $u_0\in C^2(\overline\Omega)$ verifying (\ref{i0}),
$u_{0,r}(R)=0$ and $u_0(R)>0$, property (\ref{i2}) is  in particular satisfied if we 
take a sufficiently large multiple of $u_0$ as initial data
(see Proposition~\ref{monot}(ii)).

\begin{theo}\label{theo13b}
Let $n\ge 3$. Consider problem (\ref{0}), where $u_0$ satisfies (\ref{i0}), (\ref{i2}) and $T<\infty$.
 \abs
(i)\ Let $\Omega=B_R\subset\R^n$. Then there exist $c,\eta>0$, such that
   \bas
	 U(x) \ge \frac{c}{|x|^2},
	\qquad \mbox{for all } x\in B_\eta\setminus \{0\}.
  \eas
(ii)\ Let $\Omega=\R^n$ and assume $B(u_0)\neq\R^n$. 
Then assertion (i) remains valid.
\end{theo}
Our last main result gives a more precise upper estimate which provides information on the space-time blowup behavior.
To this end, let us recall that any blowup solution of problem (\ref{0}) satisfies $\liminf_{t\to T} \, (T-t)\|u(t)\|_\infty>0$ and that
blowup is said to be type I if 
$$\limsup_{t\to T} \, (T-t)\|u(t)\|_\infty<\infty$$
 and type~II otherwise. This classification is motivated by scale invariance considerations
 and the underlying ODE. Indeed, substituting the equation for $v$ into the equation for $u$ in (\ref{0}), we obtain (in the case $\mu=0$)
  \be{NLHconv}
u_t=\Delta u+u^2-\nabla v\cdot\nabla u,
\ee
whose spatially homogeneous solutions are given by $u(t)=(T-t)^{-1}$.
Examples of type II blowup are known for $n=2$ (see \cite{HV}, \cite{RS}) and $n\ge 11$ (see \cite{MS}).
In dimensions $3\le n\le 9$, for $\Omega=\R^n$, it was proved in \cite{MS2} that blowup is type~I
provided $B(u_0)\neq\R^n$. See also \cite{GP} for a sufficient condition for type~I blowup in the case $\Omega=B_R$ with $n\ge 3$
(where $\mu=0$ and Dirichlet instead of Neumann conditions are taken for $v$ in (\ref{0})).

\begin{theo}\label{theo3}
Under the assumptions of Theorem \ref{theo13}, there exists $K>0$ such that
  \be{spacetime0}
  u(x,t) \le\Bigl({1\over u(0,t)}+K|x|^2\Bigr)^{-1},	\qquad 0<|x|\le R,\ 0<t<T.
    \ee
In particular, if blowup is type~I, then there exists $C>0$ such that
  \be{spacetime}
  u(x,t) \le C\bigl(T-t+|x|^2\bigr)^{-1},	\qquad 0<|x|\le R,\ 0<t<T.
  \ee
This is true for instance if $3\le n\le 9$ with $\Omega=\R^n$ and $u_0\in L^1(\R^n)$.
\end{theo}

\smallskip

\mysection{Remarks and discussion}
\smallskip

\hskip 3mm (i) As noted by several authors (see, e.g., \cite{biler1995}, \cite{HV}), there is a natural parallel between system~(\ref{0}) and the classical semilinear heat equation with quadratic nonlinearity
  \be{NLH}
  u_t=\Delta u+u^2,
  \ee
since the equation (\ref{NLHconv}), obtained from (\ref{0}), is precisely (\ref{NLH}) with an added convection term. 
The latter is expected to produce spreading effects
and, interestingly, Theorems~\ref{theo13} and \ref{theo13b} show that the final blowup profile is different and less singular than for (\ref{NLH}). 
Indeed, if $n\le 6$, for any radial decreasing solution of (\ref{NLH}), the blowup profile satisfies 
  \be{NLHprofile}
  U(x)\sim 16\, |x|^{-2}|\log|x||,\quad x\to 0
    \ee
(see \cite{HV92}, \cite{Vel92}, \cite{MZ98b},  \cite{Sou18pre}),
and some radial decreasing solutions with this behavior exist for all $n\ge 1$ (see \cite{BE}, \cite{Liu}, \cite{Sou18pre}).
On the other hand, for $7\le n\le 15$, beside solutions 
with profile (\ref{NLHprofile}), equation (\ref{NLH}) also possesses (radial decreasing) self-similar solutions 
of the form (\ref{Selfsimil}), leading to the homogeneous profile (\ref{Homogprofile})
(see \cite{GKS}, \cite{Lep1}, \cite{BQ89}, \cite{Lep2}).

\smallskip

\hskip 3mm (ii) In \cite{GP} and \cite{GMS}, it is shown, respectively in the case $\Omega=B_R$ 
(where  $\mu=0$ and Dirichlet instead of Neumann conditions are taken for $v$) and $\Omega=\R^n$, that
some classes of radial decreasing solutions of (\ref{0}) (satisfying suitable zero number properties) are attracted by self-similar solutions, in the sense that
$$\lim_{t\to T} (T-t)u(y\sqrt{T-t})=V(y),\quad\hbox{ uniformly for $y$ bounded}$$
(where $V$ is the profile of one of the self-similar solutions mentioned above).
 However, this convergence in the ``microscopic'' scale $y$ bounded 
does not provide information on the final blowup profile in the original variable $x$.

\smallskip

\hskip 3mm (iii) It remains an open question whether the limit $\lim_{x\to 0} |x|^2U(x)$  exists
under the assumptions of Theorems~\ref{theo13} and \ref{theo13b}.
Note that these results show that the liminf is positive and the limsup is finite.
In any case, this limit, if it exists, cannot be universal (i.e., independent of the solution, like in property (\ref{NLHprofile}) for 
equation~(\ref{NLH})).
Indeed, it is shown in \cite{GMS} that there exists a one-parameter family $u_\alpha$ of self-similar solutions of the form (\ref{Selfsimil}),
such that the limit $L(\alpha):= \lim_{y\to\infty} y^2 V_\alpha(y)=\lim_{x\to 0} |x|^2U_\alpha(x)\in (0,\infty)$ 
and is different for each $\alpha$.
\smallskip

\hskip 3mm (iv) 
The self-similar solutions in (\ref{Selfsimil}) satisfy 
$$c_1\bigl(T-t+|x|^2\bigr)^{-1}\le u(x,t) \le c_2\bigl(T-t+|x|^2\bigr)^{-1},\quad x\in\R^n,\ 0<t<T$$
(due to $\lim_{y\to\infty}y^2V(y)\in (0,\infty)$).
We expect that, for general radial nondecerasing solutions, the lower space-time estimate corresponding to (\ref{spacetime}), i.e.:
$$u(x,t) \ge c\bigl(T-t+|x|^2\bigr)^{-1},$$
should be true as $(x,t)\to (0,T)$ if blowup is type I. However we are presently unable to show this.
Note that for radially decreasing solutions of the semilinear heat equation (\ref{NLH}) with $n\le 6$,
the refined space-time blowup behavior is completely known, given by
$$ u(x,t)
= (1+o(1))\biggl[T-t+{|x|^2\over 8\min\bigl\{|\log(T-t)|,2|\log |x||\bigr\}}\biggr]^{-1},
\quad\hbox{ as $(x,t)\to (0,T)$}$$
(see \cite{HV92}, \cite{Vel92}, \cite{MZ98b}, \cite{Sou18pre}).
 This relies on powerful techniques based, among other things, on 
suitable linearization arguments, which require a variational structure that seems to be absent
for system~(\ref{0}) (even under the reduced scalar form (\ref{e1}) below).

\smallskip

\hskip 3mm (v) The upper estimates in Theorem \ref{theo13}(i) remain true for $n=2$, but they are then immediate consequences of the mass conservation property. Indeed, the latter (cf.~(\ref{mass}) below) yields
\be{masssingul}
  r^nu(r,t)\le c\int_0^r s^{n-1} u(s,t)\, ds\le c\|u(t)\|_{L^1(B_R)}=c\|u_0\|_{L^1(B_R)},\quad 0<r\le R,\ 0<t<T.
\ee
Moreover, in view of (\ref{DiracCv}), the lower estimate in Theorem~\ref{theo13b} fails for $n=2$,
and the upper estimate (\ref{13.2}) is no longer optimal at $t=T$.
 As for the space-time estimate (\ref{spacetime0}) in Theorem~\ref{theo3}, 
it also remains true for $n=2$ and may however be new in this case.

On the other hand, for $n\ge 3$, property (\ref{masssingul}) yields $u(r,t)\le Cr^{-n}$
(hence single-point blow-up for $u$, as noted in \cite[p.2140]{GP}),
but this estimate is not optimal in view of Theorem~\ref{theo13}.

\smallskip

\hskip 3mm (vi) The conclusions of Theorems \ref{theo13}-\ref{theo3} remain valid,
with same proofs, for the modified problem~(\ref{0})
(considered for instance in \cite{BHN}, \cite{GP}),
where $\Omega=B_R$, $\mu=0$ and the boundary conditions are replaced with
$\frac{\partial u}{\partial\nu}-u\frac{\partial v}{\partial\nu}=0$ and $v=0$ on $\partial\Omega$.

\smallskip

\hskip 3mm (vii) We do not know if $B(u_0)=\R^n$ can occur for radial nonincreasing solutions,
unless $u_0$ is constant. The known proof (see \cite{Chen90}) for the nonlinear heat equation (\ref{NLH}) does not seem to apply to this case.
Also, in Theorem \ref{theo13b}, the assumption $B(u_0)\neq\R^n$ is actually not necessary
since, in case $B(u_0)=\R^n$, then $U(x)\equiv\infty$ (recalling $u_t\ge 0$ and $u_r\le 0$).

\medskip

The outline of the rest of paper is as follows.
In section~3 we collect a number of preliminary results about local existence-uniqueness, 
transformed equation, blow-up criterion, and monotonicity properties.
Theorems \ref{theo13} and \ref{theo3} are then proved in Section~4, and
Theorem \ref{theo13b} in Section~5.

\bigskip

\mysection{Preliminary results}
We rewrite problem (\ref{0}) under the form
\be{0b}
        \left\{ \begin{array}{lcll}
    	\hfill u_t &=& \Delta u - \nabla \cdot (u\nabla v), 
    	& x\in\Omega, \ t>0, \\[1mm]
    	\hfill v&=&\displaystyle\int_\Omega G(x,y)(u(y,t) - \mu)\,dy,
    	& x\in\Omega, \ t>0, \\
    	\hfill \displaystyle\frac{\partial u}{\partial\nu}&=&0,
    	& x\in\pO, \ t>0, \\[2.5mm]
    	\hfill u(x,0)&=&u_0(x),
    	& x\in\Omega,
        \end{array} \right.
\ee
where $\mu$ is defined in (\ref{defmu}) and $G$ is the Neumann Green kernel of 
 $-\Delta$ in $\Omega=B_R$
(or the Newtonian potential in case $\Omega=\R^n$).
Note that problems (\ref{0}) and (\ref{0b}) are equivalent up to the addition of a constant to $v(\cdot,t)$,
and considering (\ref{0b}) instead of (\ref{0}) enables one to avoid the related uniqueness issues.
The solution of (\ref{0}) considered in Sections~1 and 2 is the solution of (\ref{0b}) given by the 
following local existence-uniqueness result.
We shall not give the proof, which is standard and follows from arguments in, e.g.,~\cite{cieslak_win}.
Here we denote by $S(t))_{t\ge 0}$ the heat semigroup on $L^\infty(\Omega)$
(with Neumann boundary conditions if $\Omega=B_R$).

\begin{prop}\label{prop12}
Let $\Omega=B_R\subset\R^n$ or $\Omega=\R^n$, with $R>0$ and $n\ge 2$,
and let $u_0$ satisfy (\ref{i0}).
There exists $\tau>0$ and a unique, classical solution $(u,v)$ of (\ref{0b}) such that
   \bas
	\left\{ \begin{array}{l}
	(u,v)\in BC^{2,1}(\overline\Omega\times (0,\tau))\times BC^{2,0}(\overline\Omega\times (0,\tau))\\[2mm]
	u-S(t)u_0\in BC(\overline\Omega\times [0,\tau)).
		\end{array} \right.
  \eas
  Moreover, $(u,v)$ can be extended to a unique maximal solution,
  whose existence time $T=T(u_0,v_0)\in (0,\infty]$ satisfies
  \be{ext}
	\mbox{either $T<\infty$ or $\displaystyle\lim_{t\to T} \|u(\cdot,t)\|_{L^\infty(\Omega)} = \infty$.}
  \ee
The couple $(u,v)$ also solves (\ref{0}) and, for each $t\in (0,T)$, the function $u(\cdot,t)$ is nonnegative and radially symmetric nonincreasing.
  Furthermore, if $\Omega=B_R$, then $u$ enjoys the mass conservation property
 \be{mass}
\|u(t)\|_{L^1(\Omega)}=\|u_0\|_{L^1(\Omega)},\quad 0<t<T,
\ee
and (\ref{mass}) remains true for $\Omega=\R^n$ if we assume in addition $u_0\in L^1(\R^n)$.
\end{prop}

\medskip

Following, e.g., \cite{BCKSV}, \cite{GP}, we shall rely on a transformed scalar equation
involving the averaged mass of $u$ over balls:

\begin{prop}\label{proptransf}
For any given solution of (\ref{0b}) provided by Proposition \ref{prop12}, let 
\be{w}
	w(r,t):=r^{-n}\int_0^r s^{n-1} u(s,t) \,ds, 
	\qquad r\in (0,R), \ t\in [0,T)
\ee
(with $R=\infty$ in case $\Omega=\R^n$). 
Then $w$ is a classical solution of
\be{e1}
w_t-w_{rr}-{n+1\over r}w_r=n(w+brw_r)(w-\tilde\mu),
	\qquad r\in (0,R), \ t\in (0,T),
\ee
where $b=1/n$ and $\tilde\mu=\mu/n$.
Moreover, for each $t\in (0,T)$, $w(\cdot,t)$ can be extended to a $C^2$ function up to $r=0$
and $w$ satisfies the boundary conditions 
\be{e1bc}
w_r(0,t)=0, \qquad w(R,t)=\tilde\mu
\ee
(dropping the second condition in (\ref{e1bc}) if $\Omega=\R^n$).
Furthermore, we have 
\be{e1wr}
w_r\le 0,	\qquad r\in (0,R), \ t\in (0,T).
\ee
\end{prop}

{\sc Proof.} 
Using $\Delta=r^{1-n}\partial_r(r^{n-1}\partial_r)$, mutiplying the second equation in (\ref{0}) with $r^{n-1}$,
integrating over $(0,r)$ and using $v_r(0,t)=0$, we get
\be{vr}
-v_r=r^{1-n}\int_0^r s^{n-1}(u-\mu)\,ds=r(w-\tilde\mu).
\ee
On the other hand, from the first equation in (\ref{0}), we obtain
$$r^{n-1}u_t-(r^{n-1}u_r)_r=-r^{n-1}u_rv_r-u(r^{n-1}v_r)_r=-(r^{n-1}uv_r)_r.$$
Integrating over $(0,r)$ and using $u_r(0,t)= v_r(0,t)=0$, we get
\be{wt}
w_t-{u_r\over r}=-{uv_r\over r}.
\ee
Differentiating (\ref{w}), we obtain $w_r=-{n\over r}w+{u\over r}$, hence
\be{uur}
u=rw_r+nw
\ee
and
\be{uur2}
u_r=rw_{rr}+(n+1)w_r.
\ee
Substituting (\ref{vr}), (\ref{uur}) and (\ref{uur2}) in (\ref{wt}), we finally obtain (\ref{e1}).

Next, from (\ref{w}), we easily see that, for each $t\in (0,T)$, $w(\cdot,t)$ can be extended to a $C^2$ function up to $r=0$.
The boundary condition (\ref{e1bc}) at $r=0$ then follows from the radial symmetry of $u$
and (\ref{uur2}).
As for the boundary condition at $r=R$, it follows from the mass conservation property~(\ref{mass}).

Finally let us prove (\ref{e1wr}). Since $u_r\le 0$ in $[0,R]\times (0,T)$ due to Proposition~\ref{prop12}, 
we have $nw(r,t)\ge r^{-n}u(r,t)\int_0^r ns^{n-1}\,ds=u(r,t)$.
Property (\ref{e1wr}) then follows from (\ref{uur}).
\qed

\begin{rem}\label{remhigherdim}
 (i) For future reference we observe that, as a consequence of Proposition~\ref{proptransf}, $w$ can be viewed as a solution of 
\be{highereq}
w_t-\tilde\Delta w=n(w+bx\cdot\tilde\nabla w)(w-\tilde\mu)\quad\hbox{ in $\tilde B\times (0,T)$},
\ee
where $\tilde B$ is the centered ball of radius $R$ in $\R^{n+2}$ (or $\tilde B=\R^{n+2}$) 
and $\tilde\Delta$, $\tilde\nabla$ are respectively the Laplacian and the spatial gradient in $n+2$ space variables.
 Also, we note that the function $z:=w-\tilde\mu\ge 0$ satisfies (cf.~(\ref{uur}):
\be{highereq2}
        \left\{ \begin{array}{lcll}
    	&&z_t-\tilde\Delta z=uz  \ge 0, 
	&\quad\hbox{ in $\tilde B\times (0,T)$},   \\[1mm]
    	&&z=0, 
	&\quad\hbox{ on $\partial\tilde B\times (0,T)$.}   \\[1mm]
        \end{array} \right.
\ee

(ii) When $\Omega=B_R$, the mass conservation property (\ref{mass}) ensures that
\be{e4a}
w(r,t)\le C(n)\|u(t)\|_1r^{-n}\le Cr^{-n}\quad\hbox{ in $(0,R]\times (0,T)$.}
\ee

\end{rem}

\medskip

For convenience, we give a short proof of blowup for large initial mass
 in the framework of radial decreasing solutions in any dimension $n\ge 2$.
A similar result (with a different proof) can be found in \cite{BHN} for a variant of system (\ref{i0})
 with Dirichlet boundary conditions on $v$.
We leave apart the case $\Omega=\R^n$ with $n=2$, for which we refer to \cite{CPZ}, \cite{BDP}
 (and we recall that more precise results are also available for $\Omega=B_R$ when $n=2$).

\begin{prop}\label{BUcond}
(i) Let $\Omega=B_R\subset\R^n$,  
with $R>0$ and $n\ge 2$,
and let $u_0$ satisfy (\ref{i0}).
There exists $c_1(n)>0$ such that if $\|u_0\|_1>c_1(n)R^{n-2}$, then the solution of (\ref{0b}) blows up in a finite time $T<\infty$.
\abs
(ii) Let $\Omega=\R^n$, $n\ge 3$, and let $u_0$ satisfy (\ref{i0}).
There exists $c_2(n)>0$ such that if $\int_{\R^n}u_0(x)e^{-|x|^2}\,dx>c(n)$,
 then the solution of (\ref{0b}) blows up in a finite time $T<\infty$.ייי
\end{prop}

\proof 
(i) Let $w$ be the corresponding solution 
of (\ref{e1})-(\ref{e1bc}) given by Proposition~\ref{proptransf}
and set $z:=w-\tilde\mu$. We have $z\ge 0$ due to (\ref{e1bc}), (\ref{e1wr}).
By Remark~\ref{remhigherdim}, $z$ satisfies
$$z_t-\tilde\Delta z
=n(z+\tilde\mu+bx\cdot \tilde\nabla z)z
=n\Bigl[z^2+\tilde\mu z+{b\over 2}x\cdot \tilde\nabla (z^2)\Bigr] 
=n\Bigl[az^2+\tilde\mu z+\displaystyle{b\over 2} \tilde\nabla\cdot (xz^2)\Bigr]$$
in $\tilde B\times (0,T)$, where $\tilde B$ is the centered ball of radius $R$ in $\R^{n+2}$ and $a=(n-2)/2n$.
Let $\phi$ be the first positive eigenfunction of $- \tilde\Delta$ in $H^1_0(\tilde B)$, normalized in $L^1(\tilde B)$,
and $\lambda_1=c(n)R^{-2}>0$ be the first eigenvalue.
Multiplying with $\phi$, integrating by parts over $\tilde B$ and using $z=0$ on $\partial{\tilde B}$, we get
$${d\over dt}\int_{\tilde B} z\phi \,dx={n-2\over 2}\int_{\tilde B} z^2\phi \,dx+(\mu-\lambda_1)\int_{\tilde B} z\phi \,dx
-\displaystyle{1\over 2}\int_{\tilde B} (x\cdot\nabla\phi)z^2 \,dx,\quad 0<t<T.$$
Note that $x\cdot \tilde\nabla\phi\le 0$ (indeed, by uniqueness, $\phi$ is radially symmetric, and it is radially decreasing
due to $(r^{n+1}\phi_r)_r=-\lambda_1 r^{n+1}\phi\le 0$ and $\phi_r(0)=0$).
Letting $y(t)=\int_{\tilde B} z\phi \,dx>0$ and using Jensen's inequality, we then obtain
$$y'(t)\ge {n-2\over 2}y^2+(\mu-\lambda_1)y,\quad 0<t<T.$$

If $n\ge 3$ and $\mu\equiv |B_R|^{-1}\|u_0\|_1\ge \lambda_1$,
then no positive solution of this differential inequality can exist globally.

If $n=2$ and $\mu>\lambda_1$, then $T=\infty$ would imply exponential growth of $z(t)$ as $t\to\infty$.
However, by (\ref{mass}) and (\ref{w}), we have
\be{e1a}
w(r,t)\le Cr^{-n}\quad\hbox{ in $(0,R]\times (0,T),$}
\ee
hence $\int_{\tilde B} z\phi \,dx\le C\int_0^R r^{-n}r^{n+1}\,dr= CR^2$: a contradiction.
Assertion~(i) follows.

\smallskip
(ii) We now have $\mu=0$, hence $z=w$. We use the above argument with 
$\tilde B$ and $\phi$ respectively replaced with $\R^{n+2}$ and $\phi=c_0e^{-|x|^2}$,
where $c_0=c_0(n)>0$ is chosen so that $\int_{\R^{n+2}}\phi\,dx=1$.
A straightforward calculation yields $\tilde\Delta\phi\ge -2(n+2)\phi$ and we arrive at
\be{kaplany}
y'(t)\ge {n-2\over 2}y^2-2(n+2)y,\quad 0<t<T
\ee
(note that all the calculations can be justified by the fast decay of $\phi$).
We thus infer $T<\infty$ whenever $n\ge 3$ and the RHS of (\ref{kaplany}) is positive at $t=0$.
In view of (\ref{w}), this can be written as $I>c(n)$, where
  \bas
I&=&\int_0^\infty r^{-n}\Bigl(\int_0^r u_0(s)s^{n-1}\,ds\Bigr)e^{-r^2} r^{n+1}\,dr \\
&=&\int_0^\infty \Bigl(\int_s^\infty re^{-r^2}\,dr\Bigr) u_0(s)s^{n-1}\,ds
={1\over 2}\int_0^\infty u_0(s) e^{-s^2}s^{n-1}\,ds.
\eas
Assertion~(ii) follows. \qed

 We finish this preliminary section with the following proposition, which recalls a standard time monotonicity property,
and also shows that large multiples of rather general initial data satisfy the assumptions of Theorem~\ref{theo13b}.

\begin{prop}\label{monot}
(i) Let $n\ge 3$. Consider problem (\ref{0b}) with $\Omega=B_R\subset\R^n$ or $\Omega=\R^n$.
Assume that $u_0$ satisfies (\ref{i0}) and (\ref{i2}).
Then the corresponding solution $w$ of (\ref{e1})-(\ref{e1bc}) satisfies $w_t\ge 0$.
\abs
(ii) Let $\Omega=B_R$ and let $\phi\in C^2(\overline\Omega)$
be radially symmetric and nonincreasing with respect to $|x|$, with $\phi(R)>0$, $\phi_r(R)=0$ and $\phi\not\equiv const$.
Then, for all $\lambda>1$ sufficiently large, $u_0=\lambda\phi$ satisfies property~(\ref{i2}).
\end{prop}

\proof
(i) Let $w_0=w(\cdot,0)$. Using (\ref{uur}), (\ref{uur2}),
we obtain, for all $r\in (0,R)$,
\be{monotinit}
\begin{array}{lcll}
        &&w_{0,rr}+\displaystyle{n+1\over r}w_{0,r}+n(w_0+brw_{0,r})(w_0-\tilde\mu) \\
&& \qquad = r^{-1}u_{0,r}+u_0\, r^{-n} 
\displaystyle\int_0^r s^{n-1}(u_0(s)-\mu)\,ds
\ge 0.
\end{array} 
\ee
In the case $\Omega= B_R$, 
the assertion then follows from a standard argument (see e.g. \cite[Proposition 52.19]{quittner_souplet} or \cite[Lem\-ma~4.4]{GP}),
that we give for completeness. 
In view of (\ref{monotinit}), 
we see that $\underline w(x,t):=w_0(|x|)$ is a subsolution of (\ref{highereq})
with $\underline w(\cdot,t)=\tilde\mu=w(\cdot,t)$ on $\partial \tilde B$.
Therefore, by the comparison principle, $w(r,\tau)\ge \underline w(r,\tau)=w_0(r)$ in $[0,R]$ for each $\tau\in (0,T)$.
Applying the comparison principle again, it follows that $w(r,\tau+t)\ge w(r,t)$ in $[0,R]$ for all $t\in (0,T- \tau)$ 
and we conclude that $w_t\ge 0$.

In the case $\Omega=\R^n$, the above argument still works provided we can apply the comparison principle,
which might cause some difficulties due to the unboundedness of the coefficient of the gradient term in
the RHS $nw^2+(x\cdot\nabla w)w$ of (\ref{highereq}).
However, owing to $w_r\le 0$, this term has a favorable sign,
since, at a possible positive maximum of the difference of two solutions $w_1, w_2$:
$$(x\cdot\nabla w_1)w_1-(x\cdot\nabla w_2)w_2=(x\cdot\nabla w_1)w_1-(x\cdot\nabla w_1)w_2=rw_{1,r}(w_1-w_2)\le 0$$
As a consequence, the required comparison principle
can be deduced from the proof of \cite[Proposition~52.6]{quittner_souplet}.

\smallskip

(ii) Set $\psi(r)=r^{-n}\int_0^r \phi(s)s^{n-1}ds$ and $\hat\mu=R^{-n}\int_0^R \phi(s)s^{n-1}ds=\lambda^{-1}\tilde\mu$. We have
\be{defz}
z(r):=\lambda^{-1}r^{1-n} \bigl[r^{n-1}u_{0,r}+u_0\int_0^r s^{n-1}(u_0(s)-\mu)\,ds\bigr]=\phi_r+\lambda r\phi(\psi-\hat\mu),
\quad 0<r\le R.
\ee
Since $\phi(r)$ is nonincreasing and nonconstant, we have %
$$\psi'(R)=-nR^{-n-1}\int_0^R \phi(s)s^{n-1}ds+R^{-1}\phi(R)
=-nR^{-n-1}\int_0^R \bigl(\phi(s)-\phi(R)\bigr)s^{n-1}ds>0.$$
By the proof of (\ref{e1wr}), $\psi(r)$ is nonincreasing. Since $\psi(R)=\hat\mu$, we deduce that, for some $C_1>0$,
\be{psipos}
\psi(r)-\hat\mu\ge C_1(R-r),\quad 0\le r\le R.
\ee
Also, since $\phi_r(0)=\phi_r(R)=0$ and $\phi\in C^2([0,R])$, we have, for some $C_2>0$,
$$\phi_r(r)\ge -C_2r(R-r),\quad 0\le r\le R.$$
Combining this with (\ref{defz}), (\ref{psipos}) and $\phi\ge \phi(R)>0$, we obtain, assuming $\lambda\ge C_2(C_1\phi(R))^{-1}$,
$$z(r)\ge \bigl[-C_2+\lambda C_1\phi(R)\bigr]r(R-r)\ge 0,\quad 0<r\le R,$$
hence (\ref{i2}).
\qed

\mysection{Upper estimates: proof of Theorems \ref{theo13} and \ref{theo3}}
\begin{rem}\label{rem3}
We shall show (see (\ref{estimw}) that any radially nonincreasing solution $w\ge 0$ of (\ref{e1}) 
which blows up only at $r=0$ satisfies 
\be{estimw2}
w(r,t) \le\Bigl({1\over w(0,t)}+Kr^2\Bigr)^{-1},\quad 0\le r<R,\ 0<t<T
\ee
 (replacing $R$ by $1$ in (\ref{estimw2}) in case $R=\infty$).
The proof works for any $b\in (0,1/2]$ (not just $b=1/n$),
or any $b>0$ if $\mu=0$.
We stress that this estimate is no longer true for $b=0$ in general,
since equation~(\ref{e1}) for $b=\mu=0$ (scaling out the factor $n$ on the RHS)
corresponds to (\ref{NLH})  in dimension $n+2$ (cf.~Section~2, Remarks (i) and (iv)).
Indeed, the proof crucially uses the fact that $b>0$ (see (\ref{e6}) below).
\end{rem}

{\sc Proof of Theorems \ref{theo13} and \ref{theo3}.}
The proof is based on nontrivial modifications of the idea in~\cite{FML}.
We shall use the maximum principle to show that a function of the form
\be{defJgen}
J=w_r+d(r)F(w)
\ee
satisfies $J\le 0$.
Specifically, $J$ will be given by
$$J=w_r+\eps rw^2$$
with $\eps>0$ small. However the somewhat tedious calculation will be more conveniently carried out by keeping the notation in (\ref{defJgen})
and choosing the functions $d$ and $F$ later.
We will denote by $C, C_1, \dots$ generic positive constants that may vary from line to line.
 In this proof we set $R:=1$ in case $\Omega=\R^n$.
\smallskip

{\bf Step 1.} {\it Single-point blowup.} 
  Note that $0\in B(u_0)$ since otherwise $u$ is uniformly bounded due to $u_r\le 0$,
contradicting (\ref{ext}).

First assume $\Omega=B_R$.
Recalling the definition (\ref{defBu0}), we have $B(u_0)=\{0\}$
as a direct consequence of (\ref{masssingul}), (\ref{e4a}) and (\ref{vr}).

Next consider the case $\Omega=\R^n$ and $u_0\in L^1(\R^n)$.
Then (\ref{masssingul}), (\ref{e4a}) and (\ref{vr}) remain valid and we conclude as before.

Finally we consider the case $\Omega=\R^n$ and $B(u_0)\neq\R^n$.
Then there exists $a>0$ such that 
$u(a,t)+|v_r(a,t)|\le C$ for all $t\in (T/2,T)$, hence $w(a,t)\le C$ by (\ref{vr}) with $\tilde\mu=0$.
Using (\ref{e1wr}), (\ref{uur}) and (\ref{w}), it follows that
$$u(r,t)\le nw(r,t)\le n(a/r)^nw(a,t)\le Cr^{-n},\quad 0<r\le a,\ T/2<t<T$$
hence, owing to $u_r\le 0$,
\be{boundu}
u(r,t)\le C\max\bigl[r^{-n},a^{-n}\bigr],\quad r>0,\ T/2<t<T.
\ee
Now combining this with (\ref{vr}) and (\ref{wt}), we get $w_t\le uw\le C\max\bigl[r^{-n},a^{-n}\bigr]w$.
Integrating in time and going back to (\ref{vr}), 
we deduce
\be{boundvr}
|v_r|=rw(r,t)\le \|w(T/2)\|_\infty\, r\exp\bigl\{CT\max\bigl[r^{-n},a^{-n}\bigr]\bigr\},\quad r>0,\ T/2<t<T.
\ee
It follows from (\ref{boundu}) and (\ref{boundvr}) that $B(u_0)=\{0\}$.

\smallskip

{\bf Step 2.} {\it Uniform negativity of $w_r$ on the parabolic boundary.} First consider the case $\Omega=B_R$. 
Setting $V=\partial_{x_1}w$, we have $V\le 0$ in $\tilde B_+=\{x\in\R^{n+2};\ |x|<R,\ x_1>0\}$, due to $w_r\le 0$.
Differentiating (\ref{highereq}), we see that $v$ solves
\be{e2}
V_t-\tilde\Delta  V=n(w+bx\cdot\nabla w) V+n(w-\tilde\mu)( (b+1)V+bx\cdot\nabla V) \quad\hbox{ in $\tilde B_+\times (0,T)$.}
\ee
It then follows from the strong maximum principle and the Hopf lemma that
$$
w_r(r, T/2)= V(r,0,\dots,0,T/2)<0\quad\hbox{ in $(0,R]$}
$$
and
 $w_{rr}(0, T/2)= V_{x_1}(0,T/2)<0$,
hence
\be{e4}
w_r(r,T/2)\le -C_1r \quad\hbox{ in $[0,R]$.}\ee
 On the other hand,
 owing to (\ref{highereq2}) and the Hopf Lemma, we have
\be{e5}
w_r(R,t)=z_r(R,t)\le -C_2\quad\hbox{ in $[T/2,T)$.}
\ee

Next consider the case $\Omega=\R^n$.
  Property (\ref{e4}) with $R=1$ 
 follows from the above argument (where now $\tilde B_+=\{x\in\R^{n+2};\ x_1>0\}$).
On the other hand,
 since $B(u_0)=\{0\}$ by Step~1, it follows from (\ref{vr}) and (\ref{uur})
that $w$ and $w_r$ are bounded on $[T/2,T)$ near $r=1$.
Property (\ref{e5}) with $R=1$ is then a consequence of the strong maximum principle applied to equation (\ref{e2}).

\smallskip
{\bf Step 3.} {\it Local parabolic inequality for $J$.}
Setting $N:=(w+brw_r)(w-\tilde\mu)$, where $b>0$, we compute
 in $(0,R)\times (0,T)$:
$$(\partial_t-\partial^2_r)(dF(w))=dF'(w)(w_t-w_{rr})-dF''(w)w_r^2-2d'F'(w)w_r-d''F(w)$$
and, differentiating (\ref{e1}),
$$(\partial_t-\partial^2_r)w_r={n+1\over r}w_{rr}-{n+1\over r^2}w_{r}+nN_r.$$
Omitting the variables $r,t,w$ when no confusion arises, it follows that
$$J_t-J_{rr}=
{n+1\over r}w_{rr}-{n+1\over r^2}w_{r}+nN_r
+dF'\Bigl({n+1\over r}w_{r}+nN\Bigr)-dF''w_r^2-2d'F'w_r-d''F.
$$
Substituting $w_r=J-dF$ and $w_{rr}=J_r-d'F-dF'w_r=J_r-dF'J+d^2FF'-d'F$, we obtain
  \bas
J_t-J_{rr}
&=&{n+1\over r}(J_r-dF'J+d^2FF'-d'F)-{n+1\over r^2}(J-dF)+nN_r \\
&&\quad+dF'\Bigl({n+1\over r}(J-dF)+nN\Bigr)-dF''(J-dF)^2-2d'F'(J-dF)-d''F.
  \eas
Setting
$${\cal P}J:=
J_t-J_{rr}-{n+1\over r}J_{r}+\Bigl({n+1\over r^2}+dF''(J-2dF)+2d'F'\Bigr)J,$$
it follows that
  \bas
{\cal P}J
&=&{n+1\over r}(d^2FF'-d'F)+{n+1\over r^2}dF+nN_r \\
&&\qquad +dF'\Bigl(-{n+1\over r}dF+nN\Bigr)
-d^3F''F^2+2dd'FF'-d''F,
  \eas
hence
\be{PJ1}
{\cal P}J={n+1\over r^2}(d-rd')F+n(N_r +dF'N)-d^3F''F^2+2dd'FF'-d''F.
\ee
On the other hand, we have
  \bas
N_r
&=&(w+brw_r)w_r+\bigl((b+1)w_r+brw_{rr}\bigr)(w-\tilde\mu)\\[1mm]
&=&\bigl(w+br(J-dF)\bigr)(J-dF)+\bigl[(b+1)(J-dF)+br(J_r-dF'J+d^2FF'-d'F)\bigr](w-\tilde\mu)\\[1mm]
&=&{\cal L}_1J+(w-bdrF)(-dF)+\bigl[-(b+1)d+br(d^2F'-d')\bigr](w-\tilde\mu)F \\[1mm]
&=&{\cal L}_1J-2dFw+b\bigl[rd^2F+(rd^2F'-rd'-d)w\bigr]F
-\tilde\mu\bigl[-(b+1)d+br(d^2F'-d')\bigr]F,
  \eas
where 
$${\cal L}_1J:=\Bigl[w+br(J-dF)-dFbr+ (b+1-bdrF')(w-\tilde\mu)\Bigr]J+br(w-\tilde\mu)J_r,$$
and
$$dF'N=dF'\bigl(w+br(J-dF)\bigr)(w-\tilde\mu)={\cal L}_2J+dF'w^2-brd^2wFF'-\tilde\mu(w-brdF)dF',$$
where ${\cal L}_2J:=dF'br(w-\tilde\mu)J$.
Therefore, setting ${\cal L}={\cal L}_1+{\cal L}_2$, we obtain
  \bas
N_r+dF'N
&=&{\cal L}J-2dFw+dF'w^2+b\bigl[rd^2F+(rd^2F'-rd'-d)w-rd^2wF'\bigr]F\\[1mm]
&&\qquad -\tilde\mu\bigl[-(b+1)dF+br(d^2FF'-d'F)+(w-brdF)dF'\bigr]\\[1mm]
&=&{\cal L}J+dw(wF'-2F)+bd\bigl[rdF-\bigl(1+\textstyle{rd'\over d}\bigr)w\bigr]F-\tilde\mu d\bigl[wF'-\bigl(1+b+b\textstyle{rd'\over d}\bigr)F\bigr].
  \eas
Combining this with (\ref{PJ1}) and setting $\tilde {\cal P}={\cal P}-n{\cal L}$, we thus have
  \bas
\tilde {\cal P}J
&=&{n+1\over r^2}(d-rd')F +ndw(wF'-2F)+nbd\bigl[rdF-\bigl(1+\textstyle{rd'\over d}\bigr)w\bigr]F\\
&&\qquad -\mu d\bigl[wF'-\bigl(1+b+b\textstyle{rd'\over d}\bigr)F\bigr] -d^3F''F^2+2dd'FF'-d''F.
  \eas
  Now choose $d=\eps r$, with $\eps>0$ to be fixed, and $F(w)=w^2$. We get
  \bas
  \tilde {\cal P}J
&=&nbd\bigl[rdF-\bigl(1+\textstyle{rd'\over d}\bigr)w\bigr]F- \mu  
d\bigl[wF'-\bigl(1+b\bigl(1+\textstyle{rd'\over d}\bigr)\bigr)F\bigr]
-d^3F''F^2+2dd'FF'\\[1mm]
&=&nb\eps rw^3\bigl[\eps r^2w-2\bigr]-\mu\eps (1-2b)rw^2-2\eps^3r^3w^4+4\eps^2rw^3,
  \eas
hence
\be{e6}
\tilde {\cal P}J=(nb-2\eps)w^3\eps r\bigl[\eps r^2w-2\bigr] -\mu\eps (1-2b)rw^2
\qquad\hbox{ in $(0,R)\times (0,T)$.}
\ee

\smallskip
{\bf Step 4.} {\it Nonlocal parabolic inequality for $J$ and conclusion.}
To relate to $J$ the term $\eps r^2w-2$ on the RHS of (\ref{e6}), we next observe that
\be{e6a}
\begin{array}{ll}	
\eps r^2w-2=2w\Bigl(\displaystyle{\eps\over 2} r^2-{1\over w}\Bigr)
&=2w\Bigl(-\displaystyle{1\over w(0,t)}+\int_0^r \Bigl(\eps s+{w_r\over w^2}(s,t)\Bigr)\,ds\Bigr) \\[4mm]
&=2w\Bigl(-\displaystyle{1\over w(0,t)}+\int_0^r w^{-2}J(s,t)\,ds\Bigr).
\end{array}
\ee
 We note that $w\ge C>0$ on $[0,R]\times [T/2,T)$. Indeed, if $\Omega=B_R$, then this is a consequence of (\ref{e1bc}), (\ref{e1wr}),
whereas if $\Omega=\R^n$ (and $R=1$), this follows from the fact that $w$ is a supersolution of the heat equation by (\ref{highereq2}).
Taking $0<\eps<nb/2$ and using (\ref{e6}), $\mu\ge 0$, and $b\le 1/2$ if $\mu>0$, it follows that
\be{e7}
\tilde {\cal P}J
\le 2\eps r(nb-2\eps)w^4\Bigl(-{1\over w(0,t)}+C_1\int_0^r J_+(s,t)\,ds\Bigr)
\ee
for some constant $C_1>0$, where $x_+=\max(x,0)$.

On the other hand, taking $\eps$ sufficiently small, it follows from (\ref{e4}), (\ref{e5}) and (\ref{e4a}) that
\be{e8}
\hbox{ $J=w_r+\eps rw^2\le 0$ on the parabolic boundary of $Q:=(0,R)\times [T_0,T)$,}
\ee
where $T_0=T/2$. Note that this remains true with $R=1$ in case $\Omega=\R^n$ by the end of Step~2
 and~(\ref{boundvr}).
Set
$$E=\bigl\{\tau\in [T_0,T);\, J\le 0 \ \hbox{on $[0,R]\times [T_0,\tau]$}\bigr\}\neq\emptyset$$
and assume for contradiction that $T_1:=\sup E<T$.
Set $T_2=(T_1+T)/2$. Then there exists $\eta>0$ such that ${1\over w(0,t)}\ge \eta$ on $[T_0,T_2]$.
Since $J\le 0$ on $[0,R]\times [T_0,T_1]$, by continuity, there exists $T_3\in (T_1,T_2)$
such that $J< \eta/(C_1R)$ on $[0,R]\times [T_0,T_3]$. It follows from (\ref{e7}) that 
\be{e9}
\tilde {\cal P}J\le 0
\quad\hbox{ in $\Sigma:=(0,R)\times (T_0,T_3]$.}
\ee
 By the definition of $\tilde {\cal P}$ in Step~3, we may write
$$\tilde {\cal P}J=J_t-J_{rr}- {n+1\over r}J_{r}+{n+1\over r^2}J - a(r,t)J
\quad\hbox{ in $\Sigma$,}$$
for some function $a\in C(\overline\Sigma)$.
Setting $\hat J:=e^{-\lambda t}J$ with $\lambda>\sup_{\Sigma}a$,
we see from (\ref{e9}) that $\hat J$ cannot attain a positive local maximum at a point $(r,t)\in(0,R)\times (T_0,T_3]$.
Indeed at such a point, we would have
$$0\le \hat J_t-\hat J_{rr}- {n+1\over r}\hat J_{r}
\le e^{-\lambda t}\bigl[\tilde {\cal P}J+ (a(r,t)-\lambda)J\bigr]<0,$$
which is impossible.
Consequently, owing to (\ref{e8}), we have $J\le 0$ in $[0,R]\times [T_0,T_3]$.
But this contradicts the definition of $T_1$. 
It follows that $T_1=T$ and we conclude from (\ref{e6a}) that 
$\eps r^2w-2\le -2w/w(0,t)$, hence
\be{estimw}
{1\over w(r,t)}\ge {\eps r^2\over 2}+ {1\over w(0,t)}
\quad\hbox{ on $[0,R]\times [T_0,T)$.}
\ee
Since $u=rw_r+nw\le nw$ and $u(0,t)=nw(0,t)$, we get
$$
{1\over u(r,t)}\ge {\eps r^2\over 2n}+ {1\over u(0,t)}
\quad\hbox{ on $[0,R]\times [T_0,T)$,}
$$
and estimate (\ref{spacetime0}) follows, hence in particular (\ref{13.1}).

Since $B(u_0)=\{0\}$, the existence of $U(x):=\lim_{t\to T}u(x,t)$ 
for all $x\in \overline\Omega\setminus\{0\}$ is an immediate consequence of interior parabolic estimates,
and we get (\ref{13.2}). 
By dominated convergence, owing to~(\ref{13.1}), we also have convergence in $L^1(B_R)$. 
The proof of Theorems \ref{theo13} and \ref{theo3} is complete.
\qed

\mysection{Lower estimate: proof of Theorem \ref{theo13b}}

\smallskip

{\sc Proof of Theorem \ref{theo13b}.}
\smallskip

{\bf Step 1.} {\it Lower estimate of $w(\cdot,T)$.}
First, similarly as in \cite[Lemma~4.5]{GP}, we have
\be{gradest}
|w_r|\leq C_1 w^{3/2}(0,t)\quad \hbox{ in $Q:=(0,R)\times (0,T)$.}
\ee
Indeed, for all $(r,t)\in Q$,
using $w_t\geq 0$, $w_r\leq 0$, $w\ge \tilde\mu\ge 0$, we obtain
$${\partial\over\partial r}\Bigl({1\over 2}w_r^2+{n\over 3}w^{3}\Bigr)
=(w_{rr}+nw^2)w_r
=\Bigl(w_t-{n+1\over r}w_r-nbrw_r(w-\tilde\mu)+\mu w\Bigr)w_r\leq 0.$$
This guarantees
$$\Bigl({1\over 2}w_r^2+{n\over 3}w^3\Bigr)(r,t)
\leq {n\over 3}w^3(0,t),$$
hence (\ref{gradest}).

We next use a modification of an argument from \cite{SouProc} (see also \cite[p.192]{quittner_souplet}).
First, since $B(u_0)=\{0\}$ by Theorem~\ref{theo13} and 
\be{wBU}
w(0,t)=\|w(t)\|_\infty\to \infty\quad\hbox{ as $t\to T$}
\ee
by (\ref{ext}) and (\ref{uur}), we may assume that 
$$w(0,t)>2w(R,t),\quad T-\delta<t<T,$$
by taking $\delta$ small enough (with $R:=1$ in the case $\Omega=\R^n$). Therefore, for all $t\in(T-\delta,T)$, there exists $r_0(t)\in (0,R)$ such that $w(r_0(t),t)=\textstyle{1\over 2} w(0,t)$. 
Note that, since $w_r<0$ in $(0,R]\times (0,T)$,
the implicit function theorem guarantees that $r_0(t)$ is unique and is a continuous function of $t$.
By (\ref{wBU}),
this implies $r_0(t)\to 0$ as $t\to T$.
In view of (\ref{gradest}), it follows that
$$-w_r\leq C_2w^{3/2},\qquad 0\leq r\leq r_0(t).$$
Integrating, we get
$$w^{-1/2}(r_0(t),t)\leq w^{-1/2}(0,t)+C_3 r_0(t)= \big(2w(r_0(t),t)\big)^{-1/2}+C_3r_0(t),$$
hence $w(r_0(t),t)\geq C_4(r_0(t))^{-2}$. Using $w_t\geq 0$, it follows that
$$w(r_0(t),T)\geq C_4(r_0(t))^{-2},\quad 0<t<T.$$
Since $r_0$ is continuous and $r_0(t)\to 0$ as $t\to T$, 
we deduce that the range $r_0((T-\delta,t))$ contains an interval of the form $(0,\eta)$, hence
$$w(r,T)\ge C_4r^{-2},\quad 0<r<\eta.$$

\smallskip
{\bf Step 2.} {\it Lower estimate of $U$.}
Going back to $U$, we have thus proved that
\be{e10}
\int_0^r s^{n-1}U(s)\,ds=r^n w(r,T) \ge C_4r^{n-2},\quad 0<r<\eta.
\ee
On the other hand, since $n\ge 3$ and $U(r)\le Cr^{-2}$ by Theorem~\ref{theo13}, it follows that
\be{e11}
\int_0^r s^{n-1}U(s)\,ds\le C\int_0^r s^{n-3}\,ds={Cr^{n-2}\over n},\quad 0<r<R.
\ee
Since $U$ is nonincreasing, combining (\ref{e10}) and (\ref{e11}), we deduce that, for each $K>0$,
$$
\begin{array}{ll}	
\displaystyle{(Kr)^nU(r)\over n}
&\ge \displaystyle\int_r^{Kr} s^{n-1}U \,ds= \int_0^{Kr} s^{n-1}U \,ds- \int_0^r s^{n-1}U\,ds \\ [3mm]
&\ge \Bigl(C_4K^{n-2}-\displaystyle{C\over n}\Bigr)r^{n-2},\qquad\qquad  0<r<\eta/K.
\end{array}
$$
Choosing $K=(2C/(nC_4))^{1/(n-2)}$ and setting $c=CK^{-n}$, we conclude that
$$U(r)\ge \bigl(nC_4K^{n-2}-C\bigr)K^{-n}r^{-2}=cr^{-2},\quad 0<r<\eta/K.$$

\qed
\vspace*{5mm}
{\bf Acknowledgement.} \quad
  The first author is partially supported by the Labex MME-DII (ANR11-LBX-0023-01).
  The second author acknowledges support of the {\em Deutsche Forschungsgemeinschaft} in the context of the project
  {\em Analysis of chemotactic cross-diffusion in complex frameworks}.
  Part of this work was done during a visit of the second author at Universit\'e Paris 13 in September 2017.

\end{document}